\documentclass[10pt]{amsart}
\input{epsf}

\usepackage{amssymb,latexsym}
\topskip  \newtheorem{theorem}{Theorem}[section]

\newtheorem{corollary}[theorem]{Corollary}

\newtheorem{lemma}[theorem]{Lemma}

\begin{document}

\pagenumbering{arabic} \pagestyle{headings}
\def\sof{\hfill\rule{2mm}{2mm}}
\def\llim{\lim_{n\rightarrow\infty}}
\def\ls{\leq}
\def\gs{\geq}
\def\LL{\frak L}
\def\qq{{\bold q}}
\def\txx{{\frac1{2\sqrt{x}}}}
\def\B{\Box}
\def\BB{\Box\hspace{-2.5pt}\Box}
\def\md{\mbox{\,{{\footnotesize mod}}\,}}

\title{A note on sum of $k$-th power of Horadam's sequence}
\maketitle

\begin{center}Toufik Mansour \footnote{Research financed by EC's
IHRP Programme, within the Research Training Network "Algebraic
Combinatorics in Europe", grant HPRN-CT-2001-00272}
\end{center}
\begin{center}
Department of Mathematics, Chalmers University of Technology,
S-412~96 G\"oteborg, Sweden

{\tt toufik@math.chalmers.se}
\end{center}
\section*{Abstract}
Let $w_{n+2}=pw_{n+1}+qw_{n}$ for $n\geq0$ with $w_0=a$ and
$w_1=b$. In this paper we find an explicit expression, in terms of
determinants,  for $\sum_{n\geq0} w_n^kx^n$ for any $k\geq1$. As a
consequence, we derive all the previously known results for this
kind of problems, as well as many new results.

\noindent{\bf Key words}: Horadam's sequence, Fibonacci numbers,
Lucas numbers, Pell numbers
\section{Introduction and the Main result}
As usual, Fibonacci numbers $F_n$, Lucas numbers $L_n$, and Pell
numbers $P_n$ are defined by the second-order linear recurrence
sequence
\begin{equation}\label{eqflp}
w_{n+2}=p w_{n+1}+q w_{n},
\end{equation}
with given $w_0=a,w_1=b$ and $n\geq0$. This sequence was
introduced, in 1965, by Horadam~\cite{H}, and it generalizes many
sequences (see~\cite{HW,HM}). Examples of such sequences are
Fibonacci numbers sequence $(F_n)_{n\geq0}$, Lucas numbers
sequence $(L_n)_{n\geq0}$, and Pell numbers sequence
$(P_n)_{n\geq0}$, when one has $p=q=b=1$, $a=0$; $p=q=b=1$, $a=2$;
and $p=2$, $q=b=1$, $a=0$; respectively. In this paper we
interested in studding the generating function for powers of
Horadam's sequence, that is,
$\mathcal{H}_k(x;a,b,p,q)=\mathcal{H}_k(x)=\sum_{n\geq0}w_n^kx^n$.

In 1962, Riordan~\cite{R} found the generating function for powers
of Fibonacci number. He proved the generating function
$\mathcal{F}_k(x)=\sum_{n\geq0}F_n^k(x)$ that satisfy the
recurrence relation
$$(1-a_kx+(-1)^kx^2)\mathcal{F}_k(x)=1+kx\sum_{j=1}^{[k/2]}(-1)^j\frac{a_{kj}}{j}\mathcal{F}_{k-2j}((-1)^ix)$$
for $k\geq1$, where $a_1=1$, $a_2=2$, $a_s=a_{s-1}+a_{s-2}$ for
$s\geq3$, and $(1-x-x^2)^{-j}=\sum a_{kj}x^{k-2j}$. Recently,
Haukkanen~\cite{Ha} studied linear combinations of Horadam's
sequences and the generating function of the ordinary product of
two of Horadam's sequences. The main result of this paper can be
formulated as follows.

Let $\Delta_k=\Delta_k(p,q)$ be the $k\times k$ matrix
$$\left(
\begin{array}{cccccc}
1-p^kx-q^kx^2&-xp^{k-1}q^1\binom{k}{1}&-xp^{k-2}q^2\binom{k}{2}&\cdots&-xp^2q^{k-2}\binom{k}{k-2}&-xpq^{k-1}\binom{k}{k-1}\\
&&&&&\\
-p^{k-1}x    &1-xp^{k-2}q^1\binom{k-1}{1}&-xp^{k-3}q^2\binom{k-1}{2}&\cdots&-xpq^{k-2}\binom{k-1}{k-2}&-xq^{k-1}\binom{k-1}{k-1}\\
&&&&&\\
-p^{k-2}x    &-xp^{k-3}q^1\binom{k-2}{1}&1-xp^{k-4}q^2\binom{k-2}{2}&\cdots&-xq^{k-2}\binom{k-2}{k-2}&0\\
&&&&&\\
-p^{k-3}x    &-xp^{k-4}q^1\binom{k-3}{1}&-xp^{k-5}q^2\binom{k-3}{2}&\cdots&0&0\\
\vdots       & \vdots                   &\vdots &  &\vdots &\vdots\\
-p^{2}x      &-xpq^1\binom{2}{1}&-xq^2\binom{2}{2}&\cdots&1&0\\
&&&&&\\
-p^{1}x      &-xq^1\binom{1}{1}&0&\cdots&0&1
\end{array}
\right),$$ and let $\delta_k=\delta_k(p,q,a,b)$ be the $k\times k$
matrix
$$\left(
\begin{array}{cccccc}
a^k+g_kx&-xp^{k-1}q^1\binom{k}{1}&-xp^{k-2}q^2\binom{k}{2}&\cdots&-xp^2q^{k-2}\binom{k}{k-2}&-xpq^{k-1}\binom{k}{k-1}\\
&&&&&\\
g_{k-1}x    &1-xp^{k-2}q^1\binom{k-1}{1}&-xp^{k-3}q^2\binom{k-1}{2}&\cdots&-xpq^{k-2}\binom{k-1}{k-2}&-xq^{k-1}\binom{k-1}{k-1}\\
&&&&&\\
g_{k-2}x    &-xp^{k-3}q^1\binom{k-2}{1}&1-xp^{k-4}q^2\binom{k-2}{2}&\cdots&-xq^{k-2}\binom{k-2}{k-2}&0\\
&&&&&\\
g_{k-3}x    &-xp^{k-4}q^1\binom{k-3}{1}&-xp^{k-5}q^2\binom{k-3}{2}&\cdots&0&0\\
\vdots       & \vdots                   &\vdots &  &\vdots &\vdots\\
g_2x &-xpq^1\binom{2}{1}&-xq^2\binom{2}{2}&\cdots&1&0\\
&&&&&\\
g_1x      &-xq^1\binom{1}{1}&0&\cdots&0&1
\end{array}
\right),$$ where $g_j=(b^j-a^jp^j)a^{k-j}$ for all
$j=1,2,\ldots,k$.
\begin{theorem}\label{thm}
The generating function $\mathcal{H}_k(x)$ is given by
$\dfrac{\det(\delta_k)}{\det(\Delta_k)}$.
\end{theorem}

The paper is organized as follows. In Section~\ref{section_2} we
give the proof of Theorem~\ref{thm} and in Section~\ref{section_3}
we give some applications for Theorem~\ref{thm}.
\section{Proofs}\label{section_2}
Let $(w_n)_{n\geq0}$ be a sequence satisfied
Relation~(\ref{eqflp}) and $k$ be any positive integer. We define a family
$\{A_{k,d}\}_{d=1}^k$ of generating functions by
\begin{equation}\label{eqaa}
A_{k,d}(x)=\sum_{n\geq0}w_{n}^{k-d}w_{n+1}^dx^{n+1}.
\end{equation}
Now we introduce two relations (Lemma~\ref{lema} and
Lemma~\ref{lemb}) between the generating functions $A_{k,d}(x)$
and $\mathcal{H}_k(x)$ that play the crucial roles in the proof of
Theorem~\ref{thm}.

\begin{lemma}\label{lema}
For any $k\geq1$,
$$(1-p^kx-q^kx^2)\mathcal{H}_k(x)
-x\sum_{j=1}^{k-1}\binom{k}{j}p^{k-j}q^jA_{k,k-j}(x)=a^k+x(b^k-a^kp^k).$$
\end{lemma}
\begin{proof}
Using the binomial theorem (see~\cite{WR}) we get
$$w_{n+2}^k=(pw_{n+1}+qw_{n})^k=p^kw_{n+1}^k+
\sum_{j=1}^{k-1}\binom{k}{j}p^{k-j}q^jw_{n+1}^{k-j}w_{n}^j
+q^kw_{n}^k.$$ Multiplying by $x^{n+2}$ and summing over all
$n\geq0$ with using Definition~(\ref{eqaa}) we have
$$\mathcal{H}_k(x)-b^kx-a^k=p^kx(\mathcal{H}_k(x)-a^k)+
x\sum_{j=1}^{k-1}\binom{k}{j}p^{k-j}q^jA_{k,k-j}(x)
+q^kx^2\mathcal{H}_k(x),$$ as requested.
\end{proof}

\begin{lemma}\label{lemb}
For any $k-1\geq d\geq1$,
$$A_{k,d}(x)-a^{k-d}b^dx=p^dx(\mathcal{H}_k(x)-a^k)+
x\sum_{j=1}^d\binom{d}{j}p^{d-j}q^jA_{k,k-j}(x).$$
\end{lemma}
\begin{proof}
Using the binomial theorem (see~\cite{WR}) we have
$$w_{n}^{k-d}w_{n+1}^d=w_n^{k-d}(pw_n+qw_{n-1})^d=
w_n^{k-d}\sum_{j=0}^d\binom{d}{j}p^{d-j}q^jw_n^{d-j}w_{n-1}^j.$$
Multiplying by $x^{n+1}$ and summing over all $n\geq0$ we get
$$A_{k,d}(x)-a^{k-d}b^dx=p^dx(\mathcal{H}_k(x)-a^k)+
x\sum_{j=1}^d\binom{d}{j}p^{k-j}q^jA_{k,k-j}(x),$$
as requested
\end{proof}


\begin{proof}(Theorem~\ref{thm})
By using the above lemmas together with definitions we get
$$\Delta_k\cdot [\mathcal{H}_k(x), A_{k,k-1}(x),A_{k,k-2}(x),\ldots,A_{k,1}(x)]^{\mbox{T}}=v_k,$$
where
$$v_k=\left[a^k+x(b^k-a^kp^k),\,(a^1b^{k-1}-p^{k-1}a^k)x,\,(a^2b^{k-2}-p^{k-2}a^k)x,\ldots,\,(a^{k-1}b^1-p^1xa^k)x
\right]^{\mbox{T}}.$$ Hence, the solution of the above equation
gives the generating function
$\mathcal{H}_k(x)=\frac{\det(\delta_k)}{\det(\Delta_k)}$, as
claimed in Theorem~\ref{thm}.
\end{proof}
\section{Applications}\label{section_3}
In this section we present some applications for
Theorem~\ref{thm}.

\subsection{Fibonacci numbers} If $a=0$ and $p=q=b=1$, then
Theorem~\ref{thm} for $k=1,2,3,4,5,6$ yields Table~1.

\begin{table}
\begin{tabular}{|l|l|l|} \hline
  $k$ & The generating function $\mathcal{H}_k(x;1,1,1,1)$ & $\mathcal{H}_k(1/100)$
  \\ \hline\hline
  1 & $\frac{x}{1-x-x^2}$ & $\frac{100}{9899}$ \\[3pt]
  2 & $\frac{x(1-x)}{(1+x)(1-3x+x^2)}$ & $\frac{9900}{979801}$ \\[3pt]
  3 & $\frac{x(1-2x-x^2)}{(1+x-x^2)(1-4x-x^2)}$ & $\frac{979900}{96940301}$ \\[3pt]
  4 & $\frac{x(1+x)(1-5x+x^2)}{(1-x)(1+3x+x^2)(1-7x+x^2)}$ & $\frac{31986700}{3161716833}$  \\[3pt]
  5 & $\frac{x(1-7x-16x^2+7x^3+x^4)}{(1-x-x^2)(1+4x-x^2)(1-11x-x^2)}$ & $\frac{9284070100}{916060399199}$  \\[3pt]
  6 & $\frac{x(1-x)(1-11x-64x^2-11x^3+x^4)}{(1+x)(1-3x+x^2)(1+7x+x^2)(1-18x+x^2)}$ & $\frac{97194791100}{9554028773189}$  \\ \hline
\end{tabular}
\begin{center} {Table 1. The generating function for the powers of Fibonacci
numbers}\end{center}
\end{table}

\subsection{Lucas numbers} If $a=2$ and $p=q=b=1$, then
Theorem~\ref{thm} for $k=1,2,3,4,5,6$ yields Table~2.

\begin{table}
\begin{tabular}{|l|l|l|} \hline
  $k$ & The generating function $\mathcal{H}_k(x;2,1,1,1)$ & $\mathcal{H}_k(1/100)$
  \\ \hline\hline
  1 & $\frac{2-x}{1-x-x^2}$ & $\frac{19900}{9899}$ \\[3pt]
  2 & $\frac{4-3x-5x^2}{(1+x)(1-3x+x^2)}$ & $\frac{3969500}{979801}$ \\[3pt]
  3 & $\frac{8-5x-36x^2+7x^3}{(1+x-x^2)(1-4x-x^2)}$ & $\frac{794640700}{96940301}$ \\[3pt]
  4 & $\frac{16-15x-180x^2+156x^3+17x^4}{(1-x)(1+3x+x^2)(1-7x+x^2)}$ & $\frac{52773853900}{3161716833}$  \\[3pt]
  5 & $\frac{32-45x-835x^2+1440x^3+745x^4-31x^5}{(1-x-x^2)(1+4x-x^2)(1-11x-x^2)}$ & $\frac{31467947446900}{916060399199}$  \\[3pt]
  6 & $\frac{64-167x-3708x^2+12323x^3+12597x^4-3188x^5-65x^6}{(1+x)(1-3x+x^2)(1+7x+x^2)(1-18x+x^2)}$ & $\frac{688573873901500}{9554028773189}$  \\ \hline
\end{tabular}
\begin{center} {Table 2. The generating function for the powers of
Lucas numbers}\end{center}
\end{table}

\subsection{Pell numbers} If $a=0$, $b=q=1$ and $p=2$, then
Theorem~\ref{thm} for $k=1,2,3,4,5,6$ yields Table~3.

\begin{table}
\begin{tabular}{|l|l|} \hline
  $k$ & The generating function $\mathcal{H}_k(x;0,1,2,1)$
  \\ \hline\hline
  1 & $\frac{x}{1-2x-x^2}$  \\[3pt]
  2 & $\frac{x(1-x)}{(1+x)(1-6x+x^2)}$  \\[3pt]
  3 & $\frac{x(1-4x-x^2)}{(1+2x-x^2)(1-14x-x^2)}$  \\[3pt]
  4 & $\frac{x(1+x)(1-14x+x^2)}{(1-x)(1+6x+x^2)(1-34x-x^2)}$  \\[3pt]
  5 & $\frac{x(1-38x-130x^2+38x^3+x^4)}{(1-2x-x^2)(1-82x-x^2)(1+14x-x^2)}$  \\[3pt]
  6 & $\frac{x(x-1)(1-104x-1210x^2-104x^3+x^4)}{(1+x)(1+34x+x^2)(1-6x+x^2)(1-198x+x^2)}$  \\ \hline
\end{tabular}
\begin{center} {Table 3. The generating function for the powers of
Pell numbers}\end{center}
\end{table}

\subsection{Chebyshev polynomials of the second kind} If $a=1$, $b=p=2t$ and $q=-1$, then
Theorem~\ref{thm} for $k=1,2,3,4,5,6$ yields Table~4.

\begin{table}
\begin{tabular}{|l|l|} \hline
  $k$ & The generating function $\mathcal{H}_k(x;1,2t,2t,-1)$
  \\ \hline\hline
  1 & $\frac{1}{1-2tx+x^2}$  \\[3pt]
  2 & $\frac{1+x}{(1-x)((1+x)^2-4xt^2)}$  \\[3pt]
  3 & $\frac{1+4tx+x^2}{(1-2tx+x^2)(1+2t(3-4t^2)x+x^2)}$\\[3pt]
  4 & $\frac{(1+x)((1-x)^2+12t^2x)}{(1-x)((1+x)^2-4t^2x)(16t^2(1-t^2)x+(1-x)^2))}$\\[3pt]
  5 & $\frac{1-6tx+2x^2+32t^3x+96t^4x^2+32t^3x^3-32t^2x^2-6x^3t+x^4}{(1+2t(3-4t^2)x+x^2)(1-2tx+x^2)(1-8t^3(4t^2-5)x-10tx+x^2)}$\\[3pt]
  6 & $\frac{(1+x)(x^4+80t^4x^3-24x^3t^2-2x^2-480t^4x^2+640t^6x^2+88t^2x^2+80t^4x-24t^2x+1)}
  {(1-x)((1+x)^2-4t^2x)((1-x)^2+16t^2(1-t^2)x)((1+x)^2-4t^2(4t^2-3)^2x)}$\\ \hline
\end{tabular}
\begin{center} {Table 4. The generating function for the powers of
Chebyshev polynomials of the second kind}\end{center}
\end{table}

More generally, if applying Theorem~\ref{thm} for $k=1,2,3,4$,
then we get the following corollary.

\begin{corollary}\
Let $k=1,2,3,4$. Then the generating function $\mathcal{H}_k(x)$
is given by $\frac{\mathcal{A}_k(x)}{\mathcal{B}_k(x)}$ where

$$\begin{array}{l}
\mathcal{A}_1=a+x(b-ap),\\
\\
\mathcal{A}_2=(a^2+xb^2)(xq-1)a^2+a^2p^2x(xq+1)-2x^2pqab,\\
\\
\mathcal{A}_3=(a^3+b^3x-a^3p^3x)(1-q^3x^2)-2xpq(a^3+b^3x)-x^2a^3p^4q+3ab^2x^2p^2q\\
\qquad+3ab^2x^3pq^3-3a^2bx^3p^2q^3+3a^2bx^2pq^2-3p^2x^2a^3q^2,\\
\\
\mathcal{A}_4=a^4+(b^4-a^4(p^4+3p^2q+q^2))x
-q(5qa^4p^4+b^4q+a^4q^3+a^4p^6+7q^2a^4p^2\\
\qquad-6qb^2a^2p^2-4b^3ap^3-4q^2ba^3p+3b^4p^2)x^2 +q^3(-8qba^3p^3-3b^4p^2+a^4q^3\\
\qquad+5qa^4p^4-6b^2a^2p^4-b^4q+a^4p^6-4q^2ba^3p+8b^3ap^3+4q^2a^4p^2+4qb^3ap)x^3\\
\qquad+q^6(ap-b)^4x^4
\end{array}$$ and
$$\begin{array}{l}
\mathcal{B}_1=1-px-x^2q,\\
\\
\mathcal{B}_2=(1+xq)(p^2x-(xq-1)^2),\\
\\
\mathcal{B}_3=(1+pqx-q^3x^2)(1-3pqx-p^3x-q^3x^2),\\
\\
\mathcal{B}_4=(1-q^2x)((1+q^2x)^2+p^2qx)((1-q^2x)^2-p^2x(p^2+4q)).
\end{array}$$
\end{corollary}

{\bf Acknowledgments}. The final version of this paper was written
while the author was visiting University of Haifa, Israel in
January 2003. He thanks the HIACS Research Center and the Caesarea
Edmond Benjamin de Rothschild Foundation Institute for
Interdisciplinary Applications of Computer Science for financial
support, and professor Alek Vainshtein for his generosity.


\begin{thebibliography}{WWW}
\bibitem[HW]{HW} G.H.~Hardy and E.M.~Wright, An introduction to the Theory of
Numbers, 4th ed. London, Oxford University Press, 1962.

\bibitem[Ha]{Ha} P.~Haukkanen, A note on Horadam's sequence,
{\em The Fibonacci Quarterly} {\bf 40:4} (2002) 358--361.

\bibitem[Ho]{H}
A.F.~Horadam, Generalization of a result of Morgado, {\em
Portugaliae Math.} {\bf 44} (1987) 131--136.

\bibitem[HM]{HM} A.F.~Horadam and J.M.~Mahon, Pell and Pell-Locas
Polynomials, {\em The Fibonacci Quarterly} {\bf 23:1} (1985) 7--20.

\bibitem[R]{R} J.~Riordan, Gereating function for powers of
Fibonacci numbers, {\em Duke Math.J.} {\bf 29} (1962) 5--12.

\bibitem[WR]{WR} E.T.~Whittaker and G.~Robinson, The Binomial Theorem, §10 in The
Calculus of Observations: A Treatise on Numerical Mathematics, 4th ed., New
York, Dover, 1967, 15--19.
\end{thebibliography}
\end{document}